\documentclass[11pt]{article}
\usepackage{setspace}
\usepackage{pstricks}
\usepackage[left=4cm,top=3cm,bottom=4cm,right=4cm]{geometry}
\usepackage{amsthm}
\usepackage{amssymb}
\theoremstyle{definition}

\newtheorem*{prf}{Proof}

\newtheorem{thm}{Theorem}
\newtheorem{defn}{Definition}
\newtheorem{lem}{Lemma}
\newtheorem{prop}{Proposition}
\newtheorem{col}{Corollary}

\title{Axiom of Neutrosophic Choice}
\author{Ahmet \c{C}evik\footnote{E-mail: a.cevik@hotmail.com} \footnote{I would like to thank the anonymous referee for sending many useful comments and suggestions.}\\ {\small Gendarmerie and Coast Guard Academy, Ankara, Turkey}}
\date{}
\begin{document}
\maketitle
\begin{abstract}
We introduce neutrosophic choice functions, the neutrosophic counterpart of the Axiom of Choice, prove some results, and discuss how it effects the foundations of mathematics in a neutrosophic setting.
\end{abstract}

\\

The Axiom of Choice is one of the axioms of Zermelo-Fraenkel Set Theory (ZFC) which is inherently used in many fundamental theorems of mathematics such as Zorn's Lemma, well-ordering theorem, even the Law of Excluded Middle (see Diaconescu \cite{Diaconescu}). The axiom is particularly of special interest for the foundations of mathematics and logic. Its non-constructive nature has led to major discussions and some objections as to whether or not we should accept this as an axiom or why we should accept it in the first place. Despite that ZFC is based on first-order classical logic, the interpretation of the axioms would slightly change if we took non-standard logics as an alternative standpoint. One of these non-standard logics is called {\em neutrosophic logic}, founded by Smarandache \cite{Smarandache1998}, \cite{Smarandache2005}, as an extension of fuzzy logic by introducing possibility of indeterminate values into attributes. Since then, neutrosophic sets and neutrosophic logic has led to an interesting theory of research with varieties of applications. For our purpose, we will be concerned with the effects of neutrosophic logic in the foundations of mathematics, but particularly in ZFC. We will primarily focus on the changes when one replaces the standard notion of choice functions with their neutrosophic counterpart. The scope of the paper will be limited to study of neutrosophic choice functions under the given interpretation and the consequences of exchanging the standard Axiom of Choice with the Axiom of Neutrosophic Choice. All set theoretical operations considered in the paper are classical.

Given a universe $U$, a neutrosophic set will contain probability degrees of {\em membership}, {\em non-membership} and {\em indeterminacy} for each $x\in U$. If $U$ is the domain of discourse and $A$ is a neutrosophic set, then for each $x\in U$, there is a probability degree that we say $x$ {\em belongs} to $A$, a degree that $x$ {\em does not belong} to $A$, and a degree that $x$ is {\em indeterminate} in $A$. So a neutrosophic set consists of triplets of probability degrees for each three attributes of each element of the domain of discourse. A neutrosophic set is defined generally as follows.

\begin{defn}
Let $U$ be a universe. A {\em neutrosophic set} $A$ is a function from $U$ to the set $\{\langle i,j,k\rangle: i,j,k\in [0,1]_{\mathbb{R}}\}$ such that $i+j+k=1$ where, for each $x\in U$, $i$ denotes the probability of $x$ {\em being a member} of $A$, $j$ denotes the probability of $x$ {\em not being a member} of $A$, and $k$ denotes the probability of $x$ {\em being indeterminate} in $A$. 
\end{defn}

Later on we will put a condition on the probability values when we introduce neutrosophic choice functions. Throughout the paper, we will assume that the reader has some familiarity with set theory. Readers may refer to Jech \cite{Jech} or Suppes \cite{Suppes} for a detailed account on axiomatic set theory. We will also use in this paper some of the notions studied in \cite{NeutrosophicComputability}. For instance, rather than treating $A$ as a function, we will treat $A$ as a set by means which will be clear in a moment. 

\begin{defn}
Let $X$ be a collection of non-empty sets. Let $f$ be a function such that $f(x)\in x$ for each $x\in X$. Then $f$ is called a {\em choice function} for $X$. 
\end{defn}

So given a collection $X$ of non-empty sets, a choice function for $X$ selects an element from each $x\in X$. We are not concerned with the construction or the explicit definition of such functions, as it is one of the axioms of set theory that such a function exists for any set. The {\em Axiom of Choice} states that every collection of non-empty sets has a choice function. Next we define the neutrosophic counterpart of choice functions. In the neutrosophic setting, we will have degrees of {\em choosing} an element, but also at the same time degrees of {\em not choosing} and leaving it {\em indeterminate}. 
% neutrosophic choice function is a natural extension of the classical choice function, i.e. it gives (1,0,0) for the chosen member of x\in X.
\begin{defn}
Let $X$ be a collection of non-empty sets and let $x\in X$. A {\em neutrosophic choice function} for $X$ is a function $f:x\rightarrow T$, where $T=\{\langle i,j,k\rangle : i,j,k\in [0,1]_{\mathbb{R}}\}$ such that $i+j+k = 1$, provided that $i\neq j\neq k$. 
\end{defn}

In the definition, every $a\in x$ is mapped to a probability distribution triplet $\langle i,j,k\rangle$ of probability values, where each index denoting respectively the probability of {\em choosing} $a$, {\em not choosing} $a$, and leaving $a$ {\em indeterminate}. The reason why we assume the sum of the probabilities of each attribute (choose, not choose, indeterminate) adds up to 1 and that are distinct from each other is because we will later get to decide only one of the attributes as the outcome of the choice function depending which attribute has the highest probability.

Let $X$ be a collection of non-empty sets and let $x\in X$. First we need to determine the conditions of a neutrosophic choice function `neutrosophically choosing' an element $a$ from $x$, leaving $a$ unchosen or leaving $a$ indeterminate as to whether or not to choose it. We will use a similar approach that was introduced in \cite{NeutrosophicComputability}. We simply pick the attribute with the greatest probability value. That is, whenever $f(a)=\langle i,j,k\rangle$, such that $i,j,k\in\mathbb{R}$ satisfying the conditions given in the definition, then we say that 

\begin{enumerate}

\item[(i)] $f$ {\em chooses} the element $a$ if $i>j$ and $i>k$;  

\item[(ii)] $f$ {\em does not choose} the element $a$ if $j>i$ and $j>k$;

\item[(iii)] $f$ {\em leaves $a$ indeterminate} if $k>i$ and $k>j$.
\end{enumerate}

Since the values in the triplet $\langle i,j,k\rangle$ are assumed to be strictly ordered, these three conditions cover all possibilities. Note that either $i$, $j$, or $k$ is the greatest element among the three.

We define $x_c$ to be the set of all elements in $x$ that are chosen by $f$, define $x_{\overline{c}}$ to be the set of all elements in $x$ that are not chosen by $f$, and finally define $x_I$ to be the set of all elements in $x$ that are left indeterminate by $f$. If the probability value of an index is zero, then the corresponding set will be empty. We can surely generalize the choice criterion to an arbitrary treshold probability value $p$ so that we say the neutrosophic choice function {\em chooses} an element if and only if the probability of choosing that element is greater than (or equal to) $p$.

Our first observation establishes the relationship between classical choice functions and neutrosophic choice functions. 
\begin{thm}
Every classical choice function is a neutrosophic choice function.
\end{thm}
\begin{prf}
Suppose that $X$ is a collection of non-empty sets. If $f$ is a classical choice function for $X$, then $f$ chooses some $a\in x$ such that $x\in X$. This implies that each $x_c$ is a singleton since $f$ chooses one element, say $a$, from each set $x$. Then, the probability of the element $a$ being chosen is non-zero and is greater than both the probability of $a$ not being chosen and the probability of $a$ being left indeterminate. Since $f$ is a classical choice function (in the sense that $f$ does not map an element to two distinct values) and since $f$ does not choose anything other than $a$ since otherwise $f$ would not be a function. We also have that for any other element $b\in x$, the probability of $b$ not being chosen is non-zero and is greater than both the probability of $b$ being chosen and the probability of $b$ being left indeterminate. This defines a neutrosophic choice function.\qed
\end{prf}

Note that the relationship is strictly one way due to the fact that a neutrosophic choice function may choose {\em more than one element} in each $x\in X$. Having this multiple choice property is one of the main differences between a classical choice function and a neutrosophic choice function. We will study the consequences of this property of choosing possibly multiple elements from each $x\in X$. Second difference is obviously the indeterminacy feature of neutrosophy which allows us to interpret basic set theoretical concepts in a non-standard way.

Now we want to examine the effects of replacing the standard Axiom of Choice with its neutrosophic counterpart. For this purpose, we shall state the Axiom of Neutrosophic Choice as follows.
\vspace{0.5cm}

\noindent{\bf Axiom of Neutrosophic Choice.} Let $X$ be a collection of non-empty sets. Then, there exists a neutrosophic choice function for $X$. 
\vspace{0.5cm}

%From now on, we shall assume this axiom in this paper and look at its consequences. 
From classical set theory, we know that given a family $\{X_i\}_{i\in I}$ of non-empty sets for some index set $I$, 
\[
\displaystyle\prod X_i \neq \emptyset.
\]

The statement above, which we shall refer to as the {\em Cartesian product statement}, is equivalent to the standard Axiom of Choice. We argue that this equivalency does not hold for the Axiom of Neutrosophic Choice.

\begin{prop}
Axiom of Neutrosophic Choice is not equivalent to the Cartesian product statement.
\end{prop}
\begin{prf}
Let $X=\{X_i\}$ be a collection of non-empty sets such that $i\in I$ for some index set $I$. It suffices to show that the existence of not every neutrosophic choice function ensures that the Cartesian product of all $X_i$'s is non-empty. A neutrosophic choice function $f$ could be defined in a way that it leaves every element of $X_i$ indeterminate.\footnote{It is worth noting that the sets $X_i$ are classical sets. The neutrosophic choice function $f$ acts as a probability distribution over the elements of $X_i$.} Note that So we suppose that the probability of each $x\in X_i$ being left indeterminate by $f$ is greater than the other two probabilities. In this case, whether the Cartesian product is empty or non-empty will be indeterminate. That is not to say that the product is empty. Nor does it mean that it is non-empty. So unlike in the classical case, Axiom of Neutrosophic Choice is not equivalent to the Cartesian product statement.\qed
\end{prf}

Despite that it looks like neutrosophic choice functions seem weaker by their own, we can use their multiple choice property to our advantage.

\begin{defn}
Let $X$ be a collection $\{X_i\}$ of non-empty sets such that $i\in I$ for some index set $I$, and let $f$ be a neutrosophic choice function for $X$. We say that $f$ has the {\em compensation property} over $X$ if for each $X_i$ that $f$ does not choose any element in $X_i$, there exists a distinct element $y\in X_j$, for $i\neq j$, such that $f$ chooses $y$ whenever $f$ also chooses some $y'\in X_j$ such that $y'\neq y$. In this case, we say that $y$ (or $y'$) {\em compensates} some of $x\in X_i$. 
\end{defn}

If $y$ compensates $x$, then we interpret this as that $y$ can be replaced by $x$. Deciding which one of $y$ and $y'$ compensates $x$ will be explained shortly (in fact we will define in the proof Theorem 2 that the one with the second highest choice probability will compensate the one with the highest choice probability). So in other words, compensation is a property which ensures that, given a collection $X$ of non-empty sets, if $x\in X$ is a set in which no element is chosen, then there exists some `compensator' element in another $y\in X$ in which the choice function chooses multiple elements.
\vspace{0.5cm}

\noindent{\bf Example. } Let $X=\{\{1,2\},\{a,b\},\{x,y,z\}\}$ be a set. Suppose that $f$ does not choose any element from $\{1,2\}$ but chooses both $y$ and $z$ from $\{x,y,z\}$. In this instance, one of $1$ or $2$ gets compensated by $y$ or $z$. So whenever we need to use some $x\in\{1,2\}$, we can instead use one of $y$ or $z$ in place of $x$.
\vspace{0.5cm}

The reason we demand to have the compensation property for neutrosophic choice functions is due to the fact that neutrosophic choice functions $f$ without this property may contradict the standard Axiom of Choice since it could be that no element is chosen by $f$  from particular elements of the given collection.

\vspace{0.5cm}

We now start discussing the consequences of replacing standard choice functions with neutrosophic choice functions. Essentially, we want all useful theorems that are true under the Axiom of Choice to be also true under the Axiom of Neutrosophic Choice as long as the neutrosophic choice function has the compensation property. There is an interesting catch though in working with compensation. We will discuss this after proving Theorem 2.

\begin{defn}
Let $\sigma$ and $\tau$ be two strings. We denote the {\em concatenation} of $\sigma$ and $\tau$ by $\sigma\tau$. If $\sigma$ is an initial segment of $\tau$, then we simply denote this by $\sigma\subset\tau$ (e.g. $011\subset0110$). If $\sigma\subset\tau$, then we say that $\tau$ is a {\em successor} of $\sigma$ (or $\sigma$ is a {\em predecessor} of $\tau$). If $\sigma\subset\tau$, then we also say that $\tau$ is an {\em extension} of $\sigma$. If $\tau$ is an infinite string, we say that $\tau$ is an {\em infinite extension} of $\sigma$. If $\tau$ is a successor of $\sigma$ and there exists no string $\upsilon$ such that $\sigma\subset\upsilon\subset\tau$, then $\tau$ is called the {\em immediate successor} of $\sigma$. If neither $\sigma\subset\tau$ nor $\tau\subset\sigma$, then we say $\sigma$ and $\tau$ are {\em incompatible}. 
\end{defn}

\begin{defn}
A {\em tree} $T$ is a set of strings such that if $\sigma\in T$ and $\tau\subset\sigma$, then $\tau\in T$. Without loss of generality we may assume that all strings are binary, i.e. if $\sigma$ is a string then $\sigma\in\{0,1\}^k$ for some $k\in\mathbb{N}$. We say that an infinite string $A$ is a {\em path} on a tree $T$ if $\sigma\in T$ for infinitely many $\sigma\subset A$.
\end{defn}

\noindent{\bf K\"{o}nig's Lemma}. Every finitely branching infinite tree has an infinite path.
\vspace{0.5cm}

We begin with showing that K\"{o}nig's Lemma holds under the Axiom of Neutrosophic Choice with the compensation property.

\begin{thm}
If for every set there exists a neutrosophic choice function with the compensation property, then K\"{o}nig's Lemma holds.
\end{thm}
\begin{prf}
Let $T$ be an infinite tree. We show that $T$ has an infinite path without assuming the standard version of the Axiom of Choice but by using the neutrosophic counterpart. At any stage $s$, let $\sigma_s$ be a string in the tree such that $\sigma_s$ has an infinite extension in $T$. That is, if $\{\sigma^k_s\}_{k\in I}$, for some index set $I$, are the immediate successors of $\sigma_s$, then $\sigma^i_s$ must have an infinite extension in $T$ for at least one $i$. Normally we would use the Axiom of Choice here to select {\em one} of the possible extensions since there may be no uniform way of choosing so. However, without assuming the standard Axiom of Choice, suppose we assume that we have a neutrosophic choice function $f$ for $T$. 
\vspace{0.5cm}

Case 1. If $f$ (neutrosophically) chooses more than one such $\sigma^i_s$, then let
\[
\tau=\textrm{max}\{p_c(\sigma^i_s)\},
\]
where $p_c(\sigma^i_s)$ denotes the probability of $\sigma^i_s$ being chosen by $f$. In other words, we let $\tau$ be the string with a maximum probability of choice. We then let $\sigma_{s+1}=\tau$. If there are more than one strings with equal probability of choice, we choose the least string among them in the lexicographical order.
\vspace{0.5cm}

Case 2. If $f$ does not get to choose any extensions of $\sigma^i_s$, then here we will rely on the compensation property of $f$. If $f$ had no compensation property, then $f$ would fail to choose any extension of $\sigma^i_s$ and we would fail to construct an infinite extension of $\sigma^i_s$. If $l$ is a level in $T$ such that $f$ chooses no string at level $l$, then we call $l$ a {\em dead level}. Before going further, we define a notion that keeps track of elements of $T$ for recording strings that compensate others in $T$. 

\begin{defn}
Let $T$ be a tree and let $\sigma\in T$. We define the {\em backward tracking} of $\sigma$, denoted by $B_\sigma$, to be the set of all $\tau\in T$ such that $\tau\subset\sigma$. Similarly, define the {\em forward tracking} of $\sigma$, denoted by $F_\sigma$, to be the set of all $\tau\in T$ such that $\sigma\subset\tau$.
\end{defn}

\begin{lem}
If $f$ is a neutrosophic choice function for $T$ such that $f$ has a compensation property, then $f$ compensates at least one string in each level $l$ of $T$ with some string at level $k$ such that $k\neq l$.
\end{lem}
\begin{prf}
Suppose $f$ is a function satisfying the hypothesis. Assume without loss of generality that there exists some level $l$ in $T$ such that $f$ chooses no string at that level. Since no string at level $l$ gets compensated by another at the same level, it must be the case that at least one string at level $l$ gets compensated by some string at level $k$ in two possible cases. First case is that if $\sigma$ is a string at level $l$ in $T$, then $f$ might have chosen some $\sigma_0$ at level $m< l$ such that $p_c(\sigma_0)<p_c(\sigma_1)$, where $\sigma_1\subset\sigma$ for some string $\sigma_1$ that $f$ has chosen along the way we constructed the tree up to the level that $\sigma$ is placed on. Second case is that $f$ has chosen no such $\sigma_0$ for any $m< l$. But since $f$ has a compensation property, there must exist some level $m>l$ such that $f$ chooses $\tau_0$ and $\tau_1$ at level $m$ such that $\tau_0\supset\sigma$ and $\tau_1\supset\sigma$ and that $\tau_0$ and $\tau_1$ are incompatible. Then, we define $\textrm{min}\{p_c(\tau_0),p_c(\tau_1)\}$ to be the string which compensates $\sigma$. This proves the lemma.
\end{prf}  

We have ensured that for each level of $T$ at which no string is chosen by $f$, there exists a string which compensates a string at that level. This suffices to construct an `extension' for every dead level using the notion of compensation.

Continuing the proof of Case 2 of Theorem 2, suppose that $l$ is a dead level and suppose $\sigma_s$ is a string at level $l$ in $T$. We need to define an extension of $\sigma_s$. If $f$ has a compensation property, then either $\sigma_s$ is compensated by some string $\tau_B\in B_\sigma$ or some $\tau_F\in F_\sigma$. To see if there exists such string $\tau_B$, we just need to look at the strings at level $m<l$ and see if $f$ has chosen $\tau_B$ whenever $f$ has chosen $\sigma_0$ at the same level $k$ as $\tau_B$ such that $\sigma_0\subset\sigma_s$. If such $\tau_B$ exists, we backtrack to $\tau_B$ and define $\sigma_{s+1}=\tau_B$. Note that we backtrack at most finitely many times since whenever we use a compensator, we mark it to prevent for future use. For a string $\sigma$ at level $l$, since there are at most $2^l$ many strings in $B_{\sigma}$ and since $f$ has the compensation property, at most finitely many strings at level $l$ will be compensated by a string in $B_{\sigma}$. This ensures that the extensions eventually grow. Suppose on the other hand that $\tau_B$ does not exist. Then by Lemma 1, $\tau_F$ must exists. In that case, we define $\sigma_{s+1}$ to be $\tau_F$. We then define $A=\bigcup_s \sigma_s$. 

We shall argue that $A$ defines an infinite path on $T$. Since $f$ is a neutrosophic choice function for $T$ with the compensation property, there exists at least one string $\sigma$ on every dead level $l$ such that $\sigma$ is compensated either by a string in $B_\sigma$ or some string in $F_\sigma$. In each case, $\sigma_{s+1}$ exists in $T$. By induction on $s$, $A$ is an infinite path on $T$. This completes the proof of the theorem.\qed
\end{prf}

One interesting question to ask is whether the constructed path is a `real' path on $T$ in the classical sense or rather an `artificial' path. The notion of compensation naturally raises questions regarding the legitimacy of the objects that are being compensated. If we want to develop a firm theory of neutrosophic logic and investigate problems revolving around the foundations of mathematics, having a sort of compensation property for neutrosophic choice functions is necessary to have in order to simulate the functionality of standard choice functions. It is also important to note that the compensation property of neutrosophic choice functions is not something that was arbitrarily made up. It is actually a natural consequence of the multiple choice property of neutrosophic choice functions. Given a family of non-empty sets $X=\{A_i\}$, if a neutrosophic function can choose possibly multiple elements from each $A_i$, when using the chosen elements during the process of any construction of an object, it is reasonable to use instead {\em many} of them whenever needed. Standard Axiom of Choice does not have the multiple choice property for the fact that standard classical functions can only map elements in the domain to a unique element in the range.

%If $X$ is a collection of non-empty sets and $f$ is a neutrosophic choice function for $X$ such that

Given a collection $X$ of non-empty sets, some neutrosophic choice functions may choose multiple elements from {\em every} $x\in X$. In particular, if $x$ is infinite then the function may be able to choose infinitely many elements from every $x\in X$. Let us call such functions, {\em neutrosophic multiple choice functions}. Such functions can be used to prove a special version of K\"{o}nig's lemma.

\begin{col}
Let $T$ be an infinitely branching infinite tree and let $f$ be a neutrosophic multiple choice function.  Then $T$ has infinitely many infinite paths.
\end{col}
\begin{prf}
In the proof of Theorem 2, at each stage $s$, given $\sigma_s$ in $T$, choose two distinct strings $\tau_0$ and $\tau_1$ in $T$ such that $\sigma_s\subset\tau_0$ and $\sigma_s\subset\tau_1$.\qed
\end{prf}

%{\bf TOPOLOJIK BASIS ICIN NOT GELECEK
%
%BAKILACAK YERLER: AXIOM OF CHOICE'UN EQUIVALENT FORMLARI: ZORN LEMMA, EVERY VECTOR SPACE HAS A BASIS, WELL-ORDERING THEOREM KANITLARI, AXIOM OF COUNTABLE CHOICE ve WEAK KONIG'S LEMMA.}

We now prove Zorn's Lemma using neutrosophic choice functions. It should be noted however that, similar as earlier, using plain neutrosophic choice functions is not sufficient to prove this. So we instead need to consider neutrosophic choice functions with the compensation property.
\vspace{0.5cm}

\noindent{\bf Zorn's Lemma}. Let $\mathcal{P}$ be a non-empty family of sets. Suppose that for each chain $\mathcal{C}$ in $\mathcal{P}$, the set $\bigcup_{C\in\mathcal{C}} C$ is in $\mathcal{P}$. Then $\mathcal{P}$ has a maximal element.
\vspace{0.5cm}

We show that, without assuming the standard Axiom of Choice, Zorn's Lemma is true under the existence of neutrosophic choice functions with the compensation property.

\begin{thm}
If there exists a neutrosophic choice function $f$ with the compensation property for every collection of non-empty sets, then Zorn's Lemma holds.
\end{thm}
\begin{prf}
Let $\mathcal{P}$ be a family of non-empty sets satisfying the hypothesis of Zorn's lemma. Suppose for a contradiction that $\mathcal{P}$ does have a maximal element. Then for every $A\in\mathcal{P}$, then set $T_A=\{Q\in\mathcal{P}: A\subset Q\}$ is non-empty. Therefore, $\{T_A\}_{A\in\mathcal{P}}$ is a family of non-empty sets. We use the axiom of neutrosophic choice (with the compensation property) to form a family of sets $\{F_A\}_{A\in\mathcal{P}}$ such that $F_A\subseteq T_A$ and $F_A$ has exactly one element for all $A\in\mathcal{P}$. How we define $\{F_A\}_{A\in\mathcal{P}}$ now needs to be explained. For every $A\in\mathcal{P}$, each element $S$ in $T_A$ will be mapped to a triplet of probability values $\langle p_c(S), p_n(S), p_I(S)\rangle$, where $p_c(S), p_n(S), p_I(S)$ denote, respectively, the probability of a neutrosophic function for choosing $S$, not choosing $S$, and leaving $S$ indeterminate. If $p_c(S)>p_n(S)>p_I(S)$ for at least one $S\in T_A$ for each $A\in\mathcal{P}$, then this basically defines a standard choice function. Suppose there exists some $A\in\mathcal{P}$ for which either $p_c(S)<p_n(S)$ or $p_c(S)<p_I(S)$ holds for all $S\in T_A$. In this case, we want to compensate $S$. We will use the notion of backward and forward tracking to keep record of every index $B\in\mathcal{P}$ such that the neutrosophic choice function chooses at least two distinct $S_1$ and $S_2$ both in $T_B$. Compare $p_c(S_1)$ and $p_c(S_2)$ and pick the least such $S_i$ with the lower probability. Then define $S_i$ to be the {\em compensator} for $S$. Then {\em mark} $S_i$ to prevent for future use. The problem is when $f$ does not choose distinct elements and that there is no element to be used as a compensator. The assumption that $f$ is a function with the compensation property however ensures the following property.
\vspace{0.5cm}

\noindent (*) For every such index $A$ for which no $S$ is chosen, there exists some index $B$ from which the neutrosophic function $f$ chooses at least two distinct sets in $T_B$ such that at least one of them is unmarked.
\vspace{0.5cm}

The property (*) ensures that for every $B\in\mathcal{P}$, $f$ either chooses some $x\in T_B$ due to that $p_c(x)>p_n(x)$ and $p_c(x)>p_I(x)$, or due to that, for $B\neq C$, some $y\in T_C$ compensates $x$. In such cases, we either define an element by neutrosophic selection or by compensation. In fact, the handling of (*) can be described computably as follows. We keep track of indices $B$ in which there exist at least two distinct sets in $B$ that are unmarked. Initially for all $B$, if $f$ chooses $x_1,x_2,\ldots,x_i\in B$ then we mark the one with the greatest choice probability. That is, we {\em mark} $x_n$ if $p_c(x_n)>p_c(x_j)$ for all $j\neq n$. We leave the rest of $x_j$ {\em unmarked}. Any unmarked element is a potential compensator. If there exists some index $B$ such that $f$ chooses no set in $T_B$, then since $f$ has the compensation property, at least one $x\in T_B$ will be compensated. We find the compensator of $x$ in a computable fashion by checking if there exists some unmarked element previously encountered. If not, we wait until the stage we encounter a new unmarked element. Such element must exist since $f$ has the compensation property. So we are ensured that the compensator for $x$ will be encountered at some further stage. So from this it follows that there is a family of sets $\{F_A\}_{A\in\mathcal{P}}$ such that $F_A\subseteq T_A$ and $F_A$ has exactly one element (either chosen by the virtue of the neutrosophic choice function directly or by compensation) for all $A\in\mathcal{P}$. For each $A\in\mathcal{P}$, let $S_A$ be the single element in $F_A$, and so $S_A\in\mathcal{P}$ and $A\subset S_A$ for all $A\in\mathcal{P}$.

Let $\mathcal{R}\subseteq\mathcal{P}$. We call $\mathcal{R}$ {\em chain-closed} if for each chain $\mathcal{C}$ in $\mathcal{R}$, the set $\bigcup_{C\in \mathcal{C}} C$ is in $\mathcal{R}$.

By hypothesis, the family $\mathcal{P}$ is chain-closed. Let $\mathcal{M}$ be the intersection of all chain-closed families in $\mathcal{P}$. Let $\mathcal{C}$ be a chain in $\mathcal{M}$. Then $\mathcal{C}$ is a chain in $\mathcal{R}$ for all chain-closed families $\mathcal{R}\subseteq \mathcal{P}$, and so $\bigcup_{C\in \mathcal{C}} C$ is in $\mathcal{R}$ for all chain-closed families $\mathcal{R}\subseteq \mathcal{P}$, and therefore $\bigcup_{C\in\mathcal{C}} C\in\mathcal{M}$. Hence, $\mathcal{M}$ is chain-closed.

Clearly, $\emptyset$ is a chain in $\mathcal{M}$ and $\bigcup_{C\in\emptyset} C=\emptyset$. Hence $\emptyset\in\mathcal{M}$ and so $\mathcal{M}$ is non-empty.

Let $A\in\mathcal{P}$ and let $\mathcal{A}=\{X\in\mathcal{P}:S_A\subseteq X\}$. Then, $\mathcal{A}\subseteq \mathcal{P}$. If $\mathcal{C}$ is a chain in $\mathcal{A}$, then $\mathcal{C}$ is a chain in $\mathcal{P}$ and so $\bigcup_{C\in\mathcal{C}} C\in \mathcal{P}$ by hypothesis. If $C\in \mathcal{C}$, then $S_A\subseteq C$, therefore $S_A\subseteq \bigcup_{C\in \mathcal{C}} C$. Hence, $\bigcup_{C\in \mathcal{C}} C\in \mathcal{A}$. So $\mathcal{A}$ is chain-closed. Therefore, $\mathcal{M}\subseteq \mathcal{A}$. But $A\not\in \mathcal{A}$ since otherwise we would have $S_A\subseteq A$ and this would contradict the fact that $A\subset S_A$. Therefore, $A\not\in\mathcal{M}$.

Since $\mathcal{M}\subseteq \mathcal{P}$ and $A\not\in\mathcal{M}$ for every $A\in\mathcal{P}$, we see that $\mathcal{M}=\emptyset$ which is a contradiction. Therefore, $\mathcal{P}$ must have a maximal element.\qed
\end{prf}

Neutrosophic choice functions without any compensation property are too weak compared to standard choice functions since it is always a possibility that a neutrosophic function leaves every element indeterminate as to whether or not to choose it. Even worse, the probability of every element not being chosen may be greater than the probability of that of the other two cases. The compensation property is a desired property of neutrosophic choice functions to make use of all features of standard choice functions. One puzzling question is to ask how compensated elements are actually included in the definition of the collection of selected objects. The choice function, after all, non-constructively defines a collection of elements being chosen. In case of neutrosophic choice functions with the compensation property, we are constructing an object by compensation. Construction by compensation allows one to define rather artificial parts of the object. To make this point clear, observe in the proof of Theorem 2 that if $l$ is a dead level, then no branch is chosen to extend the finite path that we have constructured up to stage $l$. To define an extension `in' $l$, we look for a compensator at some level $k\neq l$ to create an artificial extension in $l$. So then how is the infinite path interpreted? Certainly, the obtained collection of branches cannot be interpreted in the classical sense. The compensator replaces a branch on every dead level. Each replacement artificially creates a proceeding object in a finite sequence of strings. We shall leave this discussion to philosophers as this antinomy, we believe, can be resolved within the philosophical community of researchers in neutrosophy.

\vspace{0.5cm}

\noindent{\bf Conclusion.}

We introduced choice functions in a neutrosophic setting and investigated the effects of the neutrosopic counterpart of the Axiom of Choice. We observe that plain neutrosophic choice functions are weaker than standard choice functions. For this reason, inspired from the natural property of multiple selection of neutrosophic choice functions, we introduced the compensation property of neutrosophic choice functions to prove certain statements which are equivalent to the standard Axiom of Choice. Allowing this property has in effect the possibility of defining an interesting approach, so called {\em choice by compensation}, in constructing an object. It is our hope that we brought into attention the interplay between neutrosophic logic and the foundations of mathematics to point out some of its interesting properties and that we hope it would encourage researchers to contribute this relatively new field containing many open questions.

\end{document}